\documentclass[11pt,twoside]{article}
\usepackage[iso]{umlaute}
\usepackage{amsmath}
\usepackage{amssymb}
\usepackage{amsxtra}
\usepackage{amsthm,amsfonts,latexsym}
\usepackage{epic}
\usepackage{eepic}
\usepackage{fancybox}
\usepackage{boxedminipage}
\usepackage{calc}
\usepackage{tabularx}
\usepackage{graphicx,color}
\usepackage{picins}
\usepackage{longtable,booktabs}
\usepackage{supertabular,longtable}
\usepackage{paralist}
\accentedsymbol{\hcirc}{ {\overset{\scriptscriptstyle \circ }{ {\rm H} }}}

\newcommand{\R}{\mathbb{R}}
\newcommand{\C}{\mathbb{C}}

\newtheorem{prop}{Proposition}

\begin{document}
\parindent=0cm
{\bf 
The Lopatinski determinant of 
small shock waves may vanish
}
\par\medskip
(H. Freist\"uhler, Konstanz, and P. Szmolyan, Vienna : February 18, 2011) 
\par\bigskip
Consider a hyperbolic system
\begin{equation}
\partial_t U+\partial_{x_1}F_1(U)+\partial_{x_2}F_2(U)
=0
\label{hcl}
\end{equation}
of (at least two) conservation laws in two space variables and a corresponding 
piecewise constant Laxian shock wave
\begin{equation}
U(x,t)=
\begin{cases}
U^-,\quad x_1<0,\\
U^+,\quad x_1>0,
\end{cases}
\label{Ushock}
\end{equation}
of speed $0$.  
The Kreiss-Majda Lopatinski determinant
\begin{equation*}
\Delta(\tau,\xi)
=
\det(R_1^-(\tau,\xi),\ldots, R_{p-1}^-(\tau,\xi),
\tau[U]+i[F^\xi (U)],
R_{p+1}^+(\tau,\xi),\ldots, R_n^+(\tau,\xi))
\end{equation*}
of \eqref{Ushock} is defined on
$$
S\equiv \{(\tau,\xi)\in \C\times \R \ :\text{  Re }\!\tau \ge 0, 
|\tau|^2+\xi^2=1\},
$$
with $F^\xi\equiv\xi F_2$ and 
$$
\{R_1^-(\tau,\xi),\ldots, R_{p-1}^-(\tau,\xi)\},\quad
\{R_{p+1}^+(\tau,\xi),\ldots, R_n^+(\tau,\xi)\}
$$
continuous bases for the (extensions to $S$ of) the stable/unstable 
spaces $E^-(\tau,\xi)$, $E^+(\tau,\xi)$ 
of
$$
A(\tau,\xi)
\equiv
(\tau I+i DF^\xi(U^\mp))(DF_1(U^\mp))^{-1}.
$$
We call a simple mode $\Lambda=\Lambda(U,N,\xi)$   
of a system \eqref{hcl} {\em M\'etivier convex}\ if 
\par \medskip
{\bf (a)} $D_U\Lambda(U,N,0)\notin\hbox{ left-Im}(D(FN)(U)-\Lambda(U,N,0)I)$
(``genuine nonlinearity") and
\par\smallskip
{\bf (b)} $D^2_{\xi}\Lambda(U,N,0)>0$.
\par\medskip
It has been proved in \cite{M} 
(and reproved in \cite{FS})
that sufficiently small Laxian shock waves
associated with a M\'etivier convex mode have 
\begin{equation}
\Delta(\tau,\xi)\neq 0 \quad \text{for all } 
(\tau,\xi)\in S,
\quad
\text{if  $p=1$ or $p=n$}. 
\label{lopneq0}
\end{equation}
\par\medskip
In this note we show:
\begin{prop}
There exist symmetric constant-multiplicity hyperbolic 
systems \eqref{hcl} of conservation laws with a M\'etivier 
convex mode such that for any sufficiently small shock wave 
associated with that mode, 
\begin{equation}
\Delta(i\sigma,\xi)=0 \quad\text{for some }\sigma,\xi \in \R
\text{ with }\sigma^2+\xi^2=1.
\label{branchdeg}
\end{equation} 
\end{prop}
{\bf Proof.} Consider a system 
\begin{align}
\partial_t u+\partial_{x_1}f_1(u)+\partial_{x_2}f_2(u)&=0,\label{feq}\\
\partial_t v+\partial_{x_1}g_1(v)+\partial_{x_2}g_2(v)&=0,\label{geq}
\end{align}
where \eqref{feq} by itself is a symmetrizable hyperbolic system of 
conservation laws with modes of constant multiplicity, among which at 
least one, $\lambda_p$, is M\'etivier convex. 
Assume that for some point $u_*$ in the state space of 
$\eqref{feq}$
and for propagation direction $N_*=(1,0)^\top$, 
this mode has zero speed, 
$$
\lambda_p(u_*,N_*,0)=0.
$$
E.\ g., system $\eqref{feq}$ could be the Euler equations for 
compressible fluid flow and $\lambda_p$ the acoustic mode.

Tune \eqref{geq} by choosing $g_1,g_2:\R^2\to\R^2$ as 
\begin{equation*}
g_1(v_1,v_2)=
\begin{pmatrix} s&0\\ 0&-s\end{pmatrix}
\begin{pmatrix} v_1\\v_2\end{pmatrix},\ \ \ \      
g_2(v_1,v_2)=
\begin{pmatrix} 0&s\\ s&0\end{pmatrix}
\begin{pmatrix} v_1\\v_2\end{pmatrix}.
\end{equation*}
with some
\begin{equation}
s>\max\{|\lambda|:\lambda \text{ eigenvalue of }Df_1(u_*)\}. 
\label{supersonic}
\end{equation}
Fix a family of small-amplitude $p$-shock waves 
$(u_\epsilon^-,u_\epsilon^+)$ of \eqref{feq} with
speed $0$ and propagation direction $N_*$, perturbing from 
$u_0^\pm=u_*$.
Augmenting $u$ to $U=(u,0)$, this family trivially induces 
a family of small-amplitude shock waves 
$(U_\epsilon^-,U_\epsilon^+)$ 
of \eqref{feq}.
Since the characteristic speeds of \eqref{geq} are $-s$ and $+s$ 
and by virtue  of \eqref{supersonic}, the augmentation increases both 
the number of outgoing and that of incoming modes by $1$, on either 
side of each shock wave. 
As \eqref{feq} and \eqref{geq} are completely independent
from each other, the Lopatinski determinant 
$\Delta_\epsilon$
of the shock wave
$(U_\epsilon^-,U_\epsilon^+)$ contains a factor 
\begin{equation*}
\delta=\det(r^-(\tau,\xi),r^+(\tau,\xi))
\end{equation*}
with $r^-(\tau,\xi),r^+(\tau,\xi)$ spanning the stable/unstable spaces
of
\begin{equation*}
b(\tau,\xi)
=
\left(\tau I 
+ i\xi\begin{pmatrix} 0&s\\ s&0\end{pmatrix}\right)
\begin{pmatrix} s&0\\ 0&-s\end{pmatrix}^{-1}.
\end{equation*}
Writing 
$$
\tau=i\sigma, \text{ and }\mu=i\beta
$$ 
for the possible eigenvalue $\mu$ of $b(\tau,\xi)$, we see that 
as
\begin{equation*}
\det\left(\sigma I + \xi \begin{pmatrix} 0&s\\ s&0\end{pmatrix}
+ \beta \begin{pmatrix} s&0\\ 0&-s\end{pmatrix}\right) 
= 
-s^2\beta^2 + (\sigma^2-s^2\xi^2),
\end{equation*}
these spaces coincide at the branchpoints given by  
\begin{equation*}
\sigma\in\{\pm s(1+s^2)^{-1/2}\}
\text{ and }
\xi\in\{\pm(1+s^2)^{-1/2}\}. 
\end{equation*}
At these (four) points, $\delta$ and thus $\Delta_\epsilon$
vanish.
\par\medskip
{\bf Remark:} 
Obviously, property \eqref{branchdeg} remains true if one modifies
$f$ and $g$ by adding arbitrary smooth functions (of $u$ and $v$)  which 
are $O(|v|^2)$.
\par
\end{document}